 \documentclass[psamsfonts]{article}
 
\usepackage{amsfonts}
\usepackage{amsmath}
\usepackage{amssymb}
\usepackage{amscd}
\usepackage{amstext}
\usepackage{amsthm}

 \title{ Unitary equivalences for essential extensions of
 $C^*$-algebras
 \thanks{Research partially supported by  NSF grant DMS 0097003 .
         AMS 2000 Subject Classification Numbers:
                         Primary 46L05.
                        Key words:
                             Unitary equivalence
                                                      \protect\\}}
\author{Huaxin Lin\\
Department of Mathematics\\
University of Oregon\\
Eugene, Oregon 97403-1222}
\date{}
\begin{document}
\maketitle

\newcommand{\CA}{$C^*$-algebra}
\newcommand{\SCA}{$C^*$-subalgebra}
\newcommand{\aue}{approximate unitary equivalence}
\newcommand{\ayue}{approximately unitarily equivalent}
\newcommand{\mops}{mutually orthogonal projections}
\newcommand{\hm}{homomorphism}
\newcommand{\pisca}{purely infinite simple \CA}
\newcommand{\andeqn}{\,\,\,\,\,\, {\rm and} \,\,\,\,\,\,}
\newcommand{\QED}{\rule{1.5mm}{3mm}}
\newcommand{\morp}{contractive completely
positive linear
 map}
\newcommand{\asmorp}{asymptotic morphism}
\newcommand{\arrow}{\rightarrow}
\newcommand{\tdsum}{\widetilde{\oplus}}
\newcommand{\pa}{\|}  
\newcommand{\ep}{\varepsilon}
\newcommand{\id}{{\rm id}}
\newcommand{\aueeps}[1]{\stackrel{#1}{\sim}}
\newcommand{\aeps}[1]{\stackrel{#1}{\approx}}
\newcommand{\dt}{\delta}
\newcommand{\yu}{\fang}
\newcommand{\ca}{{\cal C}_1}
\newcommand{\Ad}{{\rm ad}}

\newtheorem{thm}{Theorem}[section]
\newtheorem{Lem}[thm]{Lemma}
\newtheorem{Prop}[thm]{Proposition}
\newtheorem{Def}[thm]{Definition}
\newtheorem{Cor}[thm]{Corollary}
\newtheorem{Ex}[thm]{Example}
\newtheorem{Pro}[thm]{Problem}
\newtheorem{Remark}[thm]{Remark}
\newtheorem{NN}[thm]{}
\renewcommand{\theequation}{e\,\arabic{section}.\arabic{equation}}
\newcommand{\N}{{\mathbb N}}

\newcommand{\Ik}{ {\cal I}^{(k)}}
\newcommand{\Iz}{{\cal I}^{(0)}}
\newcommand{\Ii}{{\cal I}^{(1)}}
\newcommand{\Ip}{{\cal I}^{(2)}}
\newcommand{\sless}{{\stackrel{\sim}{<}}}
\newcommand{\rforal}{{\rm \,\,\,for\,\,\,all\,\,\,}}

\begin{abstract}
Let $A$ be a unital separable \CA\, and $B=C\otimes {\cal K},$
where $C$ is a unital \CA. Let $\tau: A\to M(B)/B$ be a weakly
unital full essential extensions of $A$ by $B.$ We show that there
is a bijection between  a quotient group of $K_0(B)$ onto the set
of strong unitary equivalence classes of weakly unital full
essential extensions $\sigma$ such that $[\sigma]=[\tau]$ in
$KK^1(A, B).$ Consequently, when this group is zero, unitarily
equivalent full essential extensions are strongly unitarily
equivalent. When $B$ is a non-unital but $\sigma$-unital simple
\CA\, with continuous scale, we also study the problem when two
approximately unitarily equivalent essential extensions are
strongly approximately unitarily equivalent. A group is used to
compute the strongly approximate unitary equivalence classes in
the same approximate unitary equivalent class of essential
extensions.
\end{abstract}

\section{Introduction}
Let ${\cal K}$ be the \CA\, of all compact operators on the
Hilbert space $l^2$ and let $X$ be a compact metric space. Suppose
that $\tau_1, \tau_2: A\to B(l^2)/{\cal K}$ are two essential
extensions, where $A=C(X).$  The Brown-Douglas-Fillmore theorem
states that $\tau_1$ and $\tau_2$ are unitarily equivalent if and
only if $[\tau_1]=[\tau_2]$ in $KK^1(A, {\cal K}).$ Here unitarily
equivalent means that there is a partial isometry $v\in
B(l^2)/{\cal K}$ such that $v^*v=\tau_2(1),$ $vv^*=\tau_1(1)$ and
$v^*\tau_1(a)v=\tau_2(a)$ for all $a\in A.$ If both are weakly
unital, i.e., $\tau_1(1)=\tau_2(1)=1_{B(l^2)/{\cal K}},$ then
there exists a unitary $U\in B(l^2)$ such that
$\pi(U)^*\tau_1(a)\pi(U)=\tau_2(a)$ for all $a\in A,$ if
$[\tau_1]=[\tau_2],$ where $\pi: B(l^2)\to B(l^2)/{\cal K}$ is the
quotient map. In other words, unitary equivalence is same as
strong unitary equivalence in this case. This fact is particularly
important in the classification of essentially normal operators.
This is no longer true if $A$ is not a commutative \CA\, in
general. Suppose that $A={\cal O}_n$ for $n\ge 3,$ it is known
that the above is false. In fact that there are $\tau_1, \tau_2:
{\cal O}_n\to B(l^2)/{\cal K}$ for which there is a unitary $u\in
B(l^2)/{\cal K}$ such that ${\rm ad}\,u\circ \tau_1=\tau_2$ but
they are not strongly unitarily equivalent.

The usual way to get around it is to stabilize the situation. By
identifying $M_2(B(l^2)/{\cal K})$ with $B(l^2)/{\cal K},$ one
obtains $v'\in M_2(B(l^2)/{\cal K})$ such that $v+v'$ is a unitary
in $U_0(M_2(B(l^2)/{\cal K}))$ so that there is a unitary $U\in
M_2(B(l^2))$ (which is again identified with $B(l^2)$) with
$\pi(U)=v+v'.$ Hence we have $\pi(U)^*\tau_1(a)\pi(U)=\tau_2(a)$
for all $a\in A.$ Since $KK^1(A,B)$ can not detect the difference
between (weak) unitary equivalence and strong unitary equivalence,
one may choose to ignore this difference by stabilizing the
situation. However, this does not changed the fact that, for
example in the above mentioned case that $A={\cal O}_n,$ there is
no unitary $U\in B(l^2)/{\cal K}$ such that
$\pi(U)^*\tau_1(a)\pi(U)=\tau_2(a)$ for all $a\in A.$ As one sees
in the classification of essential normal operators, it is vital
sometimes to have  the strong unitary equivalence
instead of
(weak) unitary equivalence.

Let us review the case that $A=C(X)$ and $B$ is a  non-unital but
$\sigma$-unital \CA. Assume that $M(B)/B$ is a purely infinite
simple \CA\, (like the case that $B={\cal K}$). Suppose that
$\tau: C(X)\to M(B)/B$ is a weakly unital essential extension. Let
$x\in K_1(M(B)/B)$ and pick a unitary $u\in M(B)/B$ such that
$[u]=x.$ Consider $\tau_1={\rm ad}\, u\circ \tau.$ Let $\xi\in X$
be a point. Consider $D=\{\tau(f): f(\xi)=0\}.$ By a result of
Pedersen (see \cite{P}), $D^{\perp}=\{a\in M(B)/B:
bd=db=0\}\not=\{0\}.$ It is a hereditary \SCA\, of $M(B)/B.$ A
result of S. Zhang says that $M(B)/B$ has real rank zero.
Consequently there is non-zero projection $p\in D^{\perp}.$ It is
clear that $p\tau(a)=\tau(a)p$ for all $a\in C(X).$ Since $M(B)/B$
is purely infinite simple, a result of Cuntz provides a unitary
$w_1\in p(M(B)/B)p$ such that $[w_1]=[u^*].$ Now put $v=w_1\oplus
(1-p).$ Then $[v]=[u^*].$ But $v\tau(a)=\tau(a)v$ for all $a\in
A.$  So ${\rm ad}\, v \circ \tau=\tau.$ Put $z=vu.$ Then ${\rm
ad}\, z\circ \tau=\tau_1.$ But $z\in U_0(M(B)/B)$ and there is a
unitary $U\in M(B)$ so that $\pi(U)=z.$ In other words we have
just sketched the proof of the fact that, in this case, unitary
equivalence is the same as strong unitary equivalence (see Theorem
1.9 of \cite{LneII}. Unfortunately this argument can not be
applied to more general \CA s $A.$

In this short note, we show that there is  a $K$-theoretical
obstruction to prevent, in general, the unitary equivalence from
being the same as strong unitary equivalence. We will describe it
and show that when this $K$-theoretical obstruction disappears,
then the unitary equivalence is indeed the same as the strong
unitary equivalence. We find that, for a fix $z\in KK^1(A,B),$
there is a bijection between a quotient of $K_0(B)$ and strong
unitary equivalence classes of full essential extensions which are
represented by $z.$ When $A=C(X),$ this quotient of $K_0(B)$ is
zero which explains $K$-theoretically the reason why these two
unitary equivalence relations are the same. This is a direct
application of the Universal Coefficient Theorem and some of more
recent results about classification of full essential extensions.

{\bf Acknowledgement}
This work is partially supported by a grant
of National Science Foundation of U.S.A. It is initiated
during the summer 2003 when the author visiting  East China Normal
University. It is also partially supported by Zhi-Jiang Professorship
of ECNU.

\section{Preliminaries }

\begin{Def}\label{D1}
{\rm Let $A$ be a unital \CA\, and $C$ be a non-unital but
$\sigma$-unital \CA. Let $\tau_1, \tau_2: A\to M(B)/B$ be two
essential extensions. Extensions $\tau_1$ and $\tau_2$ are said to
be {\it unitarily equivalent} if there exists a partial isometry
$v\in M(B)/B$ such that $v^*v=\tau_2(1_A),$ $vv^*=\tau_1(1_A)$ and
$$
v^*\tau_1(a)v=\tau_2(a)\rforal a\in A.
$$
The extension $\tau_1$ is said to be {\it weakly unital} if
$\tau_1(1_A)=1_{M(B)/B}.$ It should be noted if $\tau_2$ is also
weakly unital and $\tau_1$ and $\tau_1$ are unitarily equivalent,
then, in the above, $v$ can be chosen to be a unitary in $M(B)/B.$

 The two essential extensions are said to be {\it strongly
unitarily equivalent} if there exists a unitary $U\in M(B)$ such
that
$$
{\rm ad}\, \pi(U)\circ \tau_1=\tau_2.
$$

It is obvious that if two essential extensions are strongly
unitarily equivalent then they are  unitarily equivalent. }
\end{Def}

\begin{Def}\label{Dfull}
{\rm Let $A$ and $D$ be two unital \CA s. A \hm\, $h: A\to D$ is said to be
{\it full} if the (closed) ideal generated by $h(a)$ is $D$ for
every nonzero $a\in A.$

Let $B$ be a non-unital but $\sigma$-unital
\CA.
An essential extension $\tau: A\to M(B)/B$ is {\it full} if $\tau$ is
a full \hm.

If $M(B)/B$ is simple, then all essential extensions are full.
If $A$ is  simple, then all weakly unital essential extensions are full.
}
\end{Def}

\begin{Def}{\rm (cf. \cite{Lnef})}\label{DP}
Let $B$ be a non-unital but $\sigma$-unital \CA.
We say $M(B)/B$ has property (P) if
for any full element $b\in M(B)/B$  there exist $x, y\in M(B)/B$
such that $xby=1.$
\end{Def}

If $M(B)/B$ is simple, then $M(B)/B$ has the property (P).
This is always the case if $B$ is purely infinite (see also Remark
\ref{PR} below). It is proved in \cite{Lnef} that, if $B=C\otimes {\cal K},$
where $C\cong C(X)$ for some finite CW complex, then $M(B)/B$
has the property (P). Other \CA s which have property (P) are discussed
in \cite{Lnef}.

Let $A$ be a separable \CA. $A$ is said to satisfy the Universal
Coefficient Theorem if for any $\sigma$-unital \CA\, $C$ one has
the following  short exact sequence
$$
0\to ext_{\mathbb Z}(K_*(A), K_{*-1}(B)) {\stackrel{\dt}\to}
KK^*(A, B){\stackrel{\Gamma}{\to}} Hom(K_*(A), K_*(B))\to 0.
$$

\CA s in the so-called ``bootstrap" class ${\cal N}$ of amenable \CA s
satisfy the Universal Coefficient Theorem (UCT) (by \cite{RS}).
When $A$ is amenable, one has $KK^1(A,B)=Ext(A,B).$

The following classification of full extensions was established in
\cite{Lnef}

\begin{thm}\label{PT}
Let $A$ be a unital separable amenable \CA\, which satisfies the
Universal Coefficient Theorem. Let $B=C\otimes {\cal K},$ where
$C$ is a unital \CA\, so that $M(B)/B$ has property {\rm (P)}.
Then two full essential extensions $\tau_1,\tau_2: A\to M(B)/B$
are unitarily equivalent if and only if
$$
[\tau_1]=[\tau_2]\,\,\,{\rm in}\,\,\, KK^1(A, B).
$$
Moreover, for any $z\in KK^1(A,B),$ there exists a full
essential extension $\tau: A\to M(B)/B$ such that $[\tau]=z.$
\end{thm}

The basic question that we consider in this short note is the
following: Suppose that $[\tau_1]=[\tau_2]$ in $KK^1(A,B)$ so that
$\tau_1$ and $\tau_2$ are unitarily equivalent. Are they strongly
unitarily equivalent? If the answer is negative in general, when
are  they strongly unitarily equivalent?

For non-stable case, the following was proved in \cite{Lneq}.

\begin{thm}\label{PT2}
Let $A$ a unital separable amenable \CA\, which satisfies the UCT
and let $B$ be a non-unital but $\sigma$-unital simple \CA\, for
which $M(B)/B$ is simple. Let $\tau_1, \tau_2: A\to M(B)/B$ be two
essential extensions. Then $\tau_1$ and $\tau_2$ are approximately
unitarily equivalent if and only if $[\tau_1]=[\tau_2]$ in $KL(A,
M(B)/B).$
\end{thm}

\begin{Remark}\label{PR}
Let $B$ be a non-unital but $\sigma$-unital simple \CA. It is
shown in \cite{Lncs2} that $M(B)/B$ is simple if and only if $B$
has continuous scale (see \cite{Lncs2}. It was also shown that when
$M(B)/B$ is simple it is purely infinite (see \cite{Lncs2}).
\end{Remark}

\begin{Def}\label{Dprod}
Let $\{A_n\}$ be a sequence of \CA s. Denote by
$l^{\infty}(\{A_n\})$ the $C^*$-product of $\{A_n\}$ and
$c_0(\{A_n\})$ the $C^*$-direct sum of $\{A_n\}.$ We will also use
$q_{\infty}(\{A_n\})$ for the quotient
$l^{\infty}(\{A_n\})/c_0(\{A_n\}).$
\end{Def}

\section{Strong  unitary equivalence}

\begin{Def}\label{DGstr}
{\rm Let $A$ be a \CA\, and $B$ be a unital \CA. Let $G\subset
U(B)$ be a normal subgroup. Suppose that $\phi_1, \phi_2: A\to B$
are \hm s. We say $\phi_1$ and $\phi_2$ are $G$-{\it strongly
unitarily equivalent} if there exists a unitary $u\in G$ such that
$$
{\rm ad}\, u\circ \phi_1(a)=\phi_2(a)\, \rforal a\in A.
$$
In the case that $G=U_0(B),$ we simply say that $\phi_1$ and
$\phi_2$ are {\it strongly unitarily equivalent}.

We say $\phi_1$ and $\phi_2$ are $G$-{strongly approximately
unitarily equivalent} if there exists a sequence of unitaries
$u_n\in G$ such that
$$
\lim_{n\to\infty}{\rm ad}\, u_n\circ \phi_2(a)=\phi_1(a)\rforal
a\in A.
$$
If $G=U_0(B),$ we simply say that $\phi_1$ and $\phi_2$ are {\it
strongly unitarily equivalent}
 }
\end{Def}

\begin{Lem}\label{L3}
Let $A$ be a separable \CA\, and $B$ be a unital \CA\, and  $h_1,
h_2: A\to B$ be two \hm s. Suppose that $G$ is a normal subgroup
of $U(B).$

{\rm (1)} Suppose that $h_2={\rm ad}\, u\circ h_1$ for some $u\in
U(B).$ Then $h_1$ and $h_2$ are $G$-strongly unitarily equivalent
if and only if there is a unitary $v\in U(B)$ such that $[u]=[v]$
in $U(B)/(U_0(B)+G)$ and  $vh_1(a)=h_1(a)v$ for all $a\in A.$

{\rm (2)}  Suppose that there exists a sequence of unitaries
$u_n\in B$ such that
$$
h_2(a)=\lim_{n\to\infty}{\rm ad}\, u_n\circ h_1(a)\rforal a\in A.
$$
Then $h_2$ and $h_1$ are $G$-strongly approximately unitarily
equivalent if and only if there exists a sequence of unitaries
$v_k\in B$ such that $[v_k]=[u_{n(k)}]$ in $U(B)/(U_0(B)+G)$ for a
subsequence $\{n(k)\} $ and
$$
\lim_{k\to\infty} \|v_kh_2(a)-h_2(a)v_k\|=0\rforal a\in A.
$$
\end{Lem}

\begin{proof}

(1)  Let $w\in G$ such that $w^*{\rm ad}\, u\circ h_1(a)w=h_1(a)$
for all $a\in A.$ Put $v=wu.$ Then $[v]=[u]\in U(B)/(U_0(B)+G)$ and
$vh_1(a)=h_1(a)v$ for all $a\in A.$

Conversely, if there is $v\in U(B)$ with $[u]=[v]$ in
$U(B)/(U_0(B)+G)$ such that $vh_1(a)=h_1(a)v$ for all $a\in A.$
Put $w=u^*v.$ Then $w\in G.$ However, we have
$$
w^*u^*h_1(a)uw=v^*h_1(a)v=h_1(a)\,\rforal a\in A.
$$

(2) Suppose that there exists a sequence of unitaries $w_n\in G$
such that
$$
\lim_{n\to\infty}{\rm ad}\, w_n\circ h_2(a)=h_1(a).
$$
Then one has two subsequences $\{m(k)\}$ and  $\{n(k)\}$ such that
$$
\lim_{n\to\infty}{\rm ad}\, w_{m(k)}u_{n(k)}\circ
h_1(a)=h_1(a)\rforal a\in A.
$$
Choose $v_k=w_{m(k)}u_{n(k)}.$ Then $[v_k]=[u_{m(k)}]$ in
$U(B)/(U_0(B)+G).$ Moreover
$$
\lim_{n\to\infty}\|v_kh_1(a)-h_1(a)v_k\|=0\rforal a\in A.
$$
Conversely, if there exists $v_k\in U(B)$ such that
$[v_k]=[u_{n(k)}]$ in $U(B)/(U_0(B)+G)$ such that
$$
\lim_{n\to\infty}\|v_kh_2(a)-h_2(a)v_k\|=0\rforal a\in A.
$$
Put $w_k=u_{n(k)}v_k^*.$ Then $w_k\in G.$

Moreover, one checks that
$$
\|{\rm ad}\, w_k\circ h_1(a)-h_2(a)\|\le \|{\rm ad}\,w_k\circ
h_1(a)- {\rm ad}\, w_k\circ h_2(a)\|+\|{\rm ad}\,w_k\circ
h_2(a)-h_2(a)\|\to 0,
$$
as $n\to\infty$ for every $a\in A.$
\end{proof}

\begin{Def}\label{Dg1}
{\rm  Let $G$ and $F$ be two groups and $u\in G$ be a
distinguished element. Define
$$
H_u(G,F)=\{x\in F: \phi(u)=x,\phi\in Hom(G,F)\}.
$$

Let $A$ and $B$ be two \CA s. Suppose that $A$ is unital. We write
\\$H_{[1_A]}(K_0(A), K_i(B))=H_1(K_0(A), K_i(B)).$
If $K_0(A)={\mathbb Z}/n{\mathbb Z}$ for some $n>1$ and
$K_i(B)={\mathbb Z},$ then $H_1(K_0(A), K_i(B))=\{0\}.$ If
$K_0(A)={\mathbb Z}$ and $[1_A]=1$ in ${\mathbb Z},$ then
$H_1(K_0(A), K_i(B))=K_i(B).$ Suppose that $K_0(A)={\mathbb
Z}/p^2{\mathbb Z}$ with $[1_A]={\bar 1}$ and $K_i(B)={\mathbb
Z}/p{\mathbb Z}\oplus {\mathbb Z}/q{\mathbb Z},$ where $(p,q)=1,$
then $H_1(K_0(A), K_i(B))={\mathbb Z}/p{\mathbb Z}.$ If $K_i(B)$
is divisible, then $H_1(K_0(A), K_i(B))=K_i(B).$ }
\end{Def}

\begin{Prop}\label{P1}
Let $A$ be a separable  \CA\,and let $C$ be a unital \CA. Let
$G\subset U(C)$ be a normal subgroup and let $h_1, h_2: A\to C$ be
two monomorphisms. Suppose that $h_1$ is unital and $h_2={\rm
ad}\, u\circ h_1.$  If $h_2$ and $h_1$ are $G$-strongly unitarily
equivalent, then $[u]\in H_1(K_0(A), K_1(C))+\bar{G},$ where
$\bar{G}$ is the image of $G$ in $K_1(C).$

\end{Prop}

\begin{proof}
Suppose that $h_1$ and $h_2$ are $G$-strongly unitarily
equivalent. Then by \ref{L3}, there is $v\in U(C)$ such that
$vh_1(a)=h_1(a)v$ and $[v]=[u]$ in $U(C)/(U_0(C)+G.$ Thus we obtain a
\hm\, $\Phi: A\otimes C(S^1)\to C$ defined by $\phi(a\otimes
f)=h_1(a)f(v)$ for $a\in A$ and $f\in C(S^1).$
Consequently, there is a
\hm\, $\phi: K_1(A\otimes C(S^1))\to K_1(C)$ such that
$[\psi(1\otimes \imath)]=[v],$ where $\imath: S^1\to S^1$ is
the identity map. Since $[(1\otimes \imath)]=[1_A]\in
K_0(A)\subset K_1(A\otimes C(S^1)),$ we obtain a \hm\, $\phi:
K_0(A)\to K_1(C)$ such that $\phi([1_A])=[u].$ This implies
$[u]\in H_1(K_0(A), K_1(C))+\bar{G}.$

\end{proof}

Suppose that $A$ satisfies the UCT and $B$ is a $\sigma$-unital
\CA. In what follows $\Gamma: KK^1(A, B)\to Hom(K_*(A),
K_{*-1}(B))$ is the surjective map given by the UCT.

\begin{Lem}\label{L1}
Let $A$ be a separable unital \CA\, in ${\cal N}$ and
$B$ be a $\sigma$-unital \CA. Let $\tau\in KK^1(A, B)$ and
$\zeta\in H_1(K_0(A), K_0(B)).$ Then, there is 
an essential extension $\tau_1: A\otimes C(S^1)\to Q(B\otimes {\cal K})$
such that
$$
[\tau_1|_{A}]=[\tau]\andeqn \Gamma([\tau_1])([1_A\otimes
\imath)]=\zeta\, \,\,\,\,{\rm (\,in} \,\,\,K_1(Q(B\otimes
{\cal K})) =K_0(B)),
$$
where $\imath(z)=z$ for $z\in S^1$ and $Q(B\otimes {\cal K})=
M(B\otimes {\cal K})/B\otimes {\cal K}.$
\end{Lem}

\begin{proof}
Let $j: A\to A\otimes C(S^1)$ be defined by $j(a)=a\otimes 1.$ If
$A\otimes C(S^1)$ is identified with $C(S^1, A),$ then $j(a)(t)=a$
for all $t\in S^1.$ By the Kunneth formula for tensor product, we
obtain $K_0(A\otimes C(S^1))=K_0(A)\oplus K_1(A)$ and
$K_1(A\otimes C(S^1))=K_1(A)\oplus K_0(A).$ Moreover, $j_{*0}: K_0(A)\to K_0(A)\oplus
K_1(A)$ may be written as $j_{*0}(x)=x\oplus 0$ and $j_{*1}:
K_1(A)\to K_1(A)\oplus K_0(A)$ may be written as
$j_{*1}(y)=y\oplus 0.$

Let $\Gamma(\tau)_i: K_i(A)\to K_{i+1}(B)$ be the map given by the UCT, $i=0,1.$
Define $\gamma_0: K_0(A\otimes C(S^1))=K_0(A)\oplus K_1(A)\to K_1(B)$  by
$\gamma_0((x,y))=\Gamma(\tau)_0(x).$

Since $\zeta\in H_1(K_0(A), K_0(B)),$ there exists
$\phi: K_0(A)\to K_0(B)$ such that $\phi([1_A])=\zeta.$
Define $\gamma_1: K_1(A\otimes C(S^1))=K_1(A)\oplus K_0(A)\to K_0(B)$
by $\gamma_1((x,y))=\Gamma(\tau)_1(x)+\phi(y).$
By the UCT, there is $\tau': A\otimes C(S^1)\to Q(B\otimes {\cal K})$
such that $\Gamma(\tau')=(\gamma_0,\gamma_1).$
Define $\psi=\tau'|_A: A\to Q(B\otimes {\cal K}).$
It follows that $\Gamma(\psi)=\Gamma(\tau'\circ j)=\Gamma(\tau).$
Let $\tau_0: A\to Q(B\otimes {\cal K})$
be an essential extension so that $[\tau_0]=[\tau]-[\psi].$
Then $\Gamma([\tau_0])=0.$ By the UCT,
$[\tau_0]\in ext_{\mathbb Z}(K_i(A), K_i(B)).$
Suppose that $[\tau_0]$ is represented by
the following two short exact sequences:
$$
0\to K_0(B)\to K_0(E)\to K_0(A)\to 0\andeqn
0\to K_1(B)\to K_1(E)\to K_1(A)\to 0.
$$
Let $G_0=K_0(E)\oplus K_1(A)$ and $G_1=K_1(E)\oplus K_0(A).$
Then $(G_0, G_1)$ gives short exact sequences
$$
0\to K_i(B)\to G_i\to K_0(A)\oplus K_1(A)\to 0
$$
It gives an element $\dt\in ext_{\mathbb Z}(K_i(A\otimes C(S^1)), K_i(B)).$
By the UCT again, there is
$\sigma: A\otimes C(S^1)\to Q(B\otimes {\cal K})$
such that $[\sigma]=\dt.$
Let $\tau_0'=\sigma|_{A}.$ Then $\tau_0'$ gives the following
short exact sequences
$$
0\to K_0(B)\to K_0(E)\to K_0(A)\to 0\andeqn
0\to K_1(B)\to K_1(E)\to K_1(A)\to 0
$$
and $[\tau_0']=[\tau_0].$
There is $\sigma': A\otimes C(S^1)\to Q(B\otimes {\cal K})$
such that $[\sigma']=-[\sigma].$
There is $\tau_1: A\otimes C(S^1)\to Q(B\otimes {\cal K})$ so that
$[\tau_1]=[\tau']-[\sigma].$
Let $\tau'':  A\otimes C(S^1)\to Q(B\otimes {\cal K})$ so that
$[\tau'']=-[\tau'].$
Then $\Gamma(\tau_1)=\Gamma(\tau').$
Moreover,
$$
[\tau]-[\tau_1|_{A}]=[\tau]+ [\tau''|_{A}]+[\sigma'|_{A}]=
[\tau]-[\psi]+ [\sigma'|_{A}]=[\tau_0]+ [\sigma'|_{A}]=0.
$$
Note that
$\Gamma([\tau_1]([1_A\otimes \imath])=\phi([1_A])=\zeta.$
\end{proof}


\begin{Lem}\label{LK1}
Let $B=C\otimes {\cal K},$ where $C$ is unital \CA\, and $M(B)/B$
has property (P). Let $u\in U(M(B)/B)$ be a unitary so that
$[u]=0$ in $K_1(M(B)/B).$ Suppose that $\phi: C(S^1)\to M(B)/B$
defined by $\phi(f)=f(u)$ for $f\in C(S^1)$ is full. Then $u\in
U_0(M(B)/B).$

Moreover, if $z\in K_1(M(B)/B),$ then there exists $u\in
U(M(B)/B)$ such that $[u]=z.$

\end{Lem}

\begin{proof}
Let $v\in U_0(M(B)/B)$ be a unitary with $sp(v)=S^1.$ Define a
monomorphism $\phi_0: C(S^1)\to M_2(M(B)/B)$ by
$\phi_0(f)=f(v)\oplus f(1)\cdot e,$ where $e=1_{M(B)/B}.$ Since
$B$ is stable, there is $z\in M_2(M(B)/B)$ such that $zz^*=e\oplus
e,$ $z^*z=e\oplus 0.$ Put $\phi_1={\rm ad}\, z\circ \phi_0.$ Put
$v_0=\phi_1(\imath),$ where $\imath$ is the identity function on
the unit circle. Since $\phi$ is full, by 2.17 of \cite{Lnef}, it is
absorbing. In particular, there exists $W\in M_2(M(B)/B)$ with
$W^*W=e\oplus 0$ and $WW^*=e\oplus e$ such that ${\rm ad}\, W\circ
\phi_1=\phi.$
 Therefore $u=W^*(u\oplus v_0)W.$
Thus, there is $W_1\in M_2(M(B)/B)$ with $W_1^*W_1=e,$
$W_1W_1^*=e\oplus e$ such that
$$
u=W_1^*(u\oplus e)W_1.
$$
 Since $B$ is stable and
$[u]=0$ in $K_1(M(B)/B),$ we may assume that ${\rm diag}(u,e)\in
U_0(M_2(M(B)/B)).$ We may write ${\rm
diag}(u,e)=\prod_{k=1}^n\exp(ia_k),$ where $a_k\in M_2(M(B)/B)$
are self-adjoint elements. It follows that
$u=W_1^*(\prod_{k=1}^n\exp(ia_k))W_1$ is connected to $W_1^*W_1=e$
by a continuous path of unitaries in $U(M(B)/B).$ So $u\in
U_0(M(B)/B).$

To see the last part of the statement, by \cite{Lnef}, there is a
full essential extension $\tau: C(S^1)\to M(B)/B$ such that
$[\tau(\imath)]=z,$ where $\imath: S^1\to S^1$ is the canonical
unitary in $C(S^1).$

\end{proof}

In the following theorem, let $z\in KK^1(A, B)$ and define
$
{\cal T}_s^1(z)
$
to be the set of strong unitary equivalence classes of
weakly unital full extensions represented by $z.$

\vspace{0.1in}

\begin{thm}\label{T1ext}
Let $A$ be a unital separable  \CA\, in ${\cal N}$
and  let $B=C\otimes {\cal K},$ where $C$ is a unital \CA\, so
that $M(B)/B$ has property (P). Let $\tau: A\to M(B)/B$ be a
weakly unital full essential extension.
Suppose that $u\in U(M(B)/B).$ Then ${\rm ad}\, u\circ \tau$ is
strongly unitarily equivalent to $\tau$ if and only if $u\in
H_1(K_0(A), K_0(B)).$

Moreover, there is a bijection
$$
\kappa:
K_0(B)/H_1(K_0(A), K_0(B))\to {\cal T}_s^1([\tau]).
$$
\end{thm}

\begin{proof}

By \ref{P1}, one needs to prove the ``if" part of the statement.

Assume that $[u]\in H_1(K_0(A), K_0(B)).$
Then, by \ref{L1}, there is an essential extension 
$\sigma: A\otimes C(S^1)\to Q(C\otimes {\cal K})$
such that $[\sigma|_{A}i]=[\tau]$ and
$[\sigma(1\otimes \imath)]=[u].$ It follows from 2.17 of \cite{Lnef}
that we may assume that $\sigma$ is full. Since $\tau$ is weakly
unital, by replacing $\sigma$ by ${\rm ad}\, w_1\circ \sigma$ for
some isometry, we may assume that $\sigma$ is also weakly unital.
By applying 2.17 of \cite{Lnef}, there is $w\in U(M(B)/B)$ such
that
$$
{\rm ad}\, w\circ \sigma|_{A}=\tau.
$$
Put $v={\rm ad}\,w\circ \sigma(1\otimes {\imath}).$ Then $[v]=[u]$
and $v\tau(a)=\tau(a)v.$ It follows from \ref{P1} that ${\rm ad}\,
u\circ \tau$ is strongly unitarily equivalent to $\tau.$

Fix $z\in KK^1(A, B)$ for which $z=[\tau].$
We define a map $\kappa: K_0(B)/H_1(K_0(A),
K_0(B))\to {\cal T}_s^1(z)$ as follows. For each $x\in K_0(B)/H_1(K_0(A), K_0(B)),$
choose $u\in U(M(B)/B)$ such that the image of $[u]$ in
$K_0(B)/H_1(K_0(A), K_0(B))$ is $x.$ Define $\kappa(x)$ to be the
strong unitary equivalence class represented by ${\rm ad}\, u\circ
\tau.$ Let $u_1, u_2\in M(B)/B.$ We have shown that ${\rm ad}\,
u_1\circ \tau$ and ${\rm ad}\, u_2\circ \tau$ are strongly
unitarily equivalent if only if $[u_1^*u_2]\in H_1(K_0(A),
K_0(B)).$ This implies that $\kappa$ is well defined and is
injective. Since for any $\sigma\in {\cal T}_s^1(z),$ there is
$u\in M(B)/B$ such that $\sigma={\rm ad}\, u\circ \tau.$ This
shows that $\kappa$ is also surjective.
\end{proof}

\begin{Cor}\label{Cext1}
Let $A$ be a separable \CA\, in ${\cal N}$
and $B=C\otimes {\cal K},$ where $C$ is unital, such that
$M(B)/B$ has property (P).
Suppose that $K_0(B)=H_1(K_0(A), K_0(B)).$ Suppose that $\tau_1,
\tau_2: A\to M(B)/B$ are two weakly unital full essential
extensions. Then they are strongly unitarily equivalent if and
only if they are unitarily equivalent.
\end{Cor}

\begin{Cor}\label{ICext2}
In \ref{Cext1}, if $K_0(A)=G\oplus {\mathbb Z}$ with $[1_A]=(0,1),$ then two
weakly unital full essential extensions are strongly unitarily
equivalent if and only if they are unitarily equivalent.
\end{Cor}

\begin{proof}
There is a \hm\, $h: K_0(A)\to {\mathbb Z}$ which maps
$[1_A]$ to $1$ in ${\mathbb Z}.$
For any $\xi\in K_0(B),$ define $\kappa: {\mathbb Z}\to K_0(B)$
by $\kappa(n)=n\xi.$ Put $\phi=\kappa\circ h.$
Thus $\xi\in H_1(K_0(A), K_0(B)).$
Therefore $K_0(B)=H_1(K_0(A), K_0(B)).$
Then \ref{Cext1} applies.

\end{proof}

\begin{Cor}\label{ICext3}
In \ref{Cext1}, if $A$ is a unital separable commutative  \CA, then two weakly
unital  full essential extensions are strongly unitarily
equivalent if and only if they are unitarily equivalent.
\end{Cor}

\begin{proof}
Since $A=C(X),$ $K_0(A)=G\oplus {\mathbb Z}$ with $[1_A]=(0, 1).$
\end{proof}

Let $A$ be a unital \CA.  If $\tau: A\to M(B)/B$ is not unital,
then there are fewer strong unitary equivalence classes of full
essential extensions $\sigma: A\to M(B)/B$ with $[\tau]=[\sigma]$
which are  not weakly unital.

\begin{Def}\label{DGp}
{\rm Let $B=C\otimes {\cal K}$ and $A$ is unital. Suppose that
$\tau: A\to M(B)/B$ is not unital and $p=\tau(1_A).$ Let
$G_p=\{z\in U(M(B)/B): z=p+v, v^*v=vv^*=1-p\}.$ Then $G_p$ is a
subgroup. Denote by $\overline{G_p}$ the image of $G_p$ in
$K_0(M(B)/B).$ Note that if $p$ is full, then there exists a
projection $e\le 1-p$ and $w\in M(B)/B$  such that $w^*w=e$ and
$ww^*=1_{M(B)/B}.$ Then $\overline{G_p}=K_1(M(B)/B).$

 Let $\sigma: A\to M(B)/B$ be another
extension. Suppose that $\tau(1_A)=p$ and $\sigma(1_A)=q.$ Note
that if $p$ is not unitarily equivalent to $q$ (in $M(B)/B$), then
$\tau$ and $\sigma$ can not be possibly strongly unitarily
equivalent.

From the following result, one also knows that if $p$ and $q$ are
unitarily equivalent and $1-p$ is full, then $\sigma$ and $\tau$
are strongly unitarily equivalent if $[\sigma]=[\tau]$ in
$KK^1(A,B).$

  }
\end{Def}

\begin{thm}\label{Tnu}
Let $A$ be a unital separable  \CA\, in ${\cal N}$
and $B=C\otimes {\cal K},$ where $C$ is unital for which
$M(B)/B$ has property (P). Let $\tau: A\to
M(B)/B$ be a full essential extension.

If $\tau(1_A)=p,$ and $p\not=1,$ then there is a bijection from
$K_0(B)/(H_1(K_0(A), K_0(B))+\overline{G_p})$ onto the set of
strong unitary equivalence classes of full essential extensions
$\sigma: A\to M(B)/B$ for which $[\sigma]=[\tau]$ and
$\sigma(1_A)$ is unitarily equivalent to $p.$

\end{thm}

\begin{proof}
Suppose that $\sigma(1_A)$ is unitarily equivalent to $p$ and
$[\sigma]=[\tau].$ Then there exists a unitary $w\in M(B)/B$ such
that $\sigma_1={\rm ad}\, w\circ \tau.$ Suppose that $[w]\in
H_1(K_0(A), K_0(B))+\overline{G_p}.$ There is $z\in U(M(B)/B)$
such that $z=p+v$ with $v^*v=vv^*=1-p$ and $[zw]\in H_1(K_0(A),
K_0(B)).$ Note that ${\rm ad} zw\circ \tau=\sigma_1.$ Thus we may
assume that $[w]\in H_1(K_0(A), K_0(B)).$ It follows from
\ref{T1ext} that $\sigma$ and $\tau$ are strongly unitarily
equivalent.

On the other hand, note since $p$ is full,
$K_1(p(M(B)/B)p)=K_1(M(B)/B)\cong K_0(B).$ Suppose that $z_1,
z_2\in K_1(p(M(B)/B)p)$ such that
$\overline{z_1}\not=\overline{z_2}$ in $K_0(B)/(H_1(K_0(A),
K_0(B))+\overline{G_p}).$ It follows from \ref{LK1} that there are
unitaries $v_1, v_2\in p(M(B)/B)p$ such that $[v_1]=z_1$ and
$[v_2]=z_2.$ Consider extensions ${\rm ad}\, v_i\circ \tau: A\to
p(M(B)/B)p$ for $i=1,2.$  It follows from \ref{P1} that they are not
$G_p$-strongly unitarily equivalent as unital \hm s to
$p(M(B)/B)p.$ It follows that they are not strong unitarily as \hm
s to $M(B)/B.$
\end{proof}






\begin{Prop}\label{Phomot}
Let $A$ be a unital separable  \CA\, in ${\cal N}$ and
$B=C\otimes {\cal K},$ where $C$ is unital such that $M(B)/B$ has
property (P).  Suppose that $K_1(M(B)/B)=H_1(K_0(A), K_1(M(B)/B).$
Suppose that  $\tau_1,\tau_2: A\to M(B)/B$ are two weakly unital
full essential extensions. Then $\tau_1$ and $\tau_2$ are homotopic
if and only if $[\tau_1]=[\tau_2]$ in $KK^1(A, B).$
\end{Prop}

\begin{proof}
If $\tau_1$ and $\tau_2$ are homotopic, then $[\tau_1]=[\tau_2]$
in $KK^1(A, B).$

Conversely  if $\tau_1$ and $\tau_2,$ by \ref{T1ext}, they are
strongly unitarily equivalent. There exists a unitary $u\in
U_0(M(B)/B)$ such that ${\rm ad}\, u\circ \tau_1=\tau_2.$ Let
$\{u_t: t\in [0,1]\}$ be a continuous path of unitaries in
$M(B)/B$ such that $u_0=u$ and $u_1=1_{M(B)/B}.$ Define $\Sigma:
A\to C([0,1], M(B)/B)$ by $\Sigma(a)(t)={\rm ad}\, u_t\circ
\tau_1$ for $a\in A.$  Then $\Sigma(a)(0)=\tau_2$ and
$\Sigma(a)(1)=\tau_1.$
\end{proof}

\section{Approximate unitary equivalences}

\begin{Def}\label{DapHa}
Let $\{x_n\}$ be a sequence of elements in $K_i(B).$ Denote by
$H_1^{ap}(K_0(A), K_i(B))$ the subset of those sequences $\{x_n\}$
of $K_i(B)$ such that there exists an increasing sequence of
finitely generated subgroups $G_n\subset K_0(A)$ with $[1_A]\in
G_n$ and group \hm s $h_n: G_n\to K_0(B)$ such that
$h_n([1_A])=x_n.$ It forms a subgroup of $\prod K_i(B).$

Suppose that $\Pi_i: \prod_{n\in N} K_i(B)\to \prod_{n\in N} K_i(B)/\oplus_{n\in N}
K_i(B),$ $i=0,1.$
Suppose that $\xi=\Pi_i(\{x_n\})$ and $\xi\in
H_1(K_0(A), \prod_{n\in N} K_i(B)/\oplus_{n\in N} K_i(B)).$
Then $\{x_n\}\in H_1^{ap}(K_0(A), K_i(B)).$

If $K_i(B)=H_1(K_0(A), K_i(B)),$ then
$H_1^{ap}(K_0(A), K_i(B))=\prod_{n\in N} K_i(B).$
\end{Def}

Recall that, for a unitary $u$ in  a unital \CA\, $A,$
$$
{\rm cel}(u)=\inf\{\sum_{k=1}^n\|h_k\|: u=\prod_{k=1}^n{\rm exp}(ih_k),\,n\in \N, h_k\in A_{s.a}\}
$$
and ${\rm cel}(A)=\sup_{u\in U(A)}{\rm cel}(u)$ (see \cite{Pc}).

Let $r: \N\to\N$ be a map. Unital \CA\, $A$ is said to have
$K_1$-$r$-cancellation if $u\oplus 1_{M_{r(n)}}$ and $v\oplus
1_{M_{r(n)}}$ are in the same path connected component of
$U(M_{n+{r(n)}})$ for any $u, v\in M_n(A)$ with $[u]=[v]$ in
$K_1(A).$

\vspace{0.1in}

\begin{Prop}\label{IIIP1}
Let $A$ be a unital separable amenable \CA, and
$C$ be a unital \CA\, and $G\subset U(B)$ be
a normal subgroup.
Suppose that ${\rm cel}(B)\le L$ for some $L>\pi$ and
$B$ has $K_1$-$r$-cancellation (for some $r: \N\to \N$).
Let $h, \phi: A\to B$ be two unital monomorphisms.
Suppose that there exists
a sequence of unitaries $u_n\in U(B)$ such that
$$
\lim_{n\to\infty}{\rm u_n}\circ h(a)=\phi(a)\rforal a\in A.
$$
Suppose also that  $h$ and $\phi$ are $G$-strongly approximately unitarily
equivalent. Then there exists a subsequence $\{n(k)\}$
such that $[u_{n(k)}]\in H_1^{ap}(K_0(A), K_1(B)) + \prod {\bar G},$
where ${\bar G}$ is the image of $G$ in $K_1(B).$

If $K_0(A)$ is finitely generated, then one can require
that $[u_{n(k)}]\in H_1(K_0(A), K_1(B)) + {\bar G}.$

\end{Prop}

\begin{proof}
There exists a sequence of unitaries
    $v_k\in B$ such that $[v_k]=[u_{n(k)}]$ in $U(B)/(U_0(B)+G)$ for a
    subsequence $\{n(k)\} $ and
    $$
    \lim_{k\to\infty} \|v_kh_2(a)-h_2(a)v_k\|=0\rforal a\in A.
    $$
Define a map $\Phi: A\otimes C(S^1)\to l^{\infty}(B)$ by
$\Phi(a\otimes f)=\{h_2(a)f(v_k)\}$ for $f\in C(S^1).$
Let $\Pi: l^{\infty}(B)\to l^{\infty}(B)/c_0(B)$ be the quotient map.
Then $\Psi=\Pi\circ \Phi: A\otimes C(S^1)\to l^{\infty}(B)/c_0(B)$
is a \hm.
Since ${\rm cel}(B)\le L$ and
$B$ has $K_1$-$r$-cancellation, it follows from Proposition 2.1 (3) of \cite{GL} that
$$
K_1(l^{\infty}(B))\subset \prod K_1(B)\andeqn
K_1(l^{\infty}(B)/c_0(B))\subset \prod K_1(B)/\oplus K_1(B).
$$
Let $\xi=\Pi_{*1}([\{v_k\}]).$
Then one obtains a \hm\, $\gamma: K_0(A)\to \prod K_1(B)/\oplus K_1(B)$
induced by $\Psi$ which maps $[1_A]$ to $\xi.$
Let $F\subset \gamma(K_0(A)).$ It is a countable
abelian group. Write $F=\cup_nF_n,$ where $F_n\subset F_{n+1}$ and
each $F_n$ is finitely generated. It is easy to see
that for every finitely generated subgroup $F_n,$ there is a \hm\,
$f_n: F_n\to \prod K_1(B)$ such that $\Pi\circ f_n={\rm id}_{F_n}.$
Let $p_n: \prod K_1(B)\to K_1(B)$ be the projection on the
$n$ coordinate.
For each $n,$ there is $m(n)$ such that
$p_k\circ f_n\circ \gamma([1_A])=[v_k]$ if $k\ge m(n).$ We may assume
that $m(n+1)>m(n).$
Let  $G_n\subset K_0(A)$  be a finitely generated subgroup
such that $\gamma(G_n)=F_n$ and $G_n\subset G_{n+1}.$
We may also assume that $[1_A]\in G_n.$
Define $\phi_n=p_{m(n)}\circ f_n\circ \gamma.$ Then $\phi_n: G_n\to K_1(B)$
is a \hm\, and $\phi_n([1_A])=[v_{m(n)}].$
Thus $\{[v_{m(n)}]\}\in H_1^{ap}(K_0(A), K_1(B)).$
It follows that $\{[u_{n(m(k))}]\}\in H_1^{ap}(K_0(A), K_1(B))+\prod {\bar G}.$

\end{proof}

\begin{thm}\label{IIITapp=}
Let $A$ be a unital separable  \CA\, in ${\cal N}.$ Let
$B=C\otimes {\cal K},$ where $C$ is unital such that $M(B)/B$ has
property (P). Suppose that $K_0(B)=H_1(K_0(A), K_0(B)).$  Then two
weakly unital full essential extensions of $A$ by $B$ are strongly
approximately  unitarily equivalent if and only if they are
approximate unitarily equivalent.
\end{thm}

\begin{proof}
Since ${\rm cel}(M(B)/B)\le {\rm cel}(M(B)),$ by 5.1.15 of
\cite{Lnb},  ${\rm cel}(M(B)/B)\le 6\pi.$ By Lemma \ref{LK1},
$M(B)/B$ also has $K_1$-$r$-cancellation (with $r:\N\to \N$ is
the identity map). Suppose that $\tau_1:
A\to M(B)/B$ is a weakly unital full essential extension. Let
$\{u_n\}\subset M(B)/B$ be a sequence of unitaries such that
$$
\tau_2(a)=\lim_{n\to\infty}{\rm ad}\, u_n\circ \tau_1(a)\rforal
a\in A.
$$
For each $n,$ since $[u_n]\in H_1(K_0(A), K_0(B)),$
as in the proof of \ref{T1ext}, one obtains a unitary $w_n\in U(M(B)/B)$
such that $[w_n]=[u_n]$ and $w_n\tau_1(a)=\tau_1(a)w_n.$
Now let $v_n=w_n^*u_n.$ Note $v_n\in U_0(M(B)/B)$ (by \ref{LK1}). Then
$$
\tau_2(a)=\lim_{n\to\infty}{\rm ad}\, v_n\circ \tau_1(a)\rforal a\in A.
$$
So $\tau_1$ and $\tau_2$ are strongly
approximately unitarily equivalent.
\end{proof}

\begin{thm}\label{IITappf}
Let $A$ be a unital separable  \CA\, in ${\cal N}.$ Let
$B=C\otimes {\cal K},$ where $C$ is unital such that $M(B)/B$ has
property (P). Suppose that
$\tau: A\to M(B)/B$ is a weakly unital full essential extension.
 Then there is an injection from
$K_0(B)/H_1(K_0(A), K_0(B))$ to the set of strongly approximately
unitarily equivalent classes of weakly unital full essential
extensions $\sigma: A\to M(B)/B$ for which $[\sigma]=[\tau]$ in
$KL^1(A, B).$
\end{thm}

\begin{proof}
Let $x_1,x_2\in K_0(B)/H_1(K_0(A),K_0(B)).$ Suppose that $u_i\in
U(M(B)/B)$ such that $\overline{[u_i]}=x_i$ in $K_0(B)/H_1(K_0(A),
K_0(B)),$ $i=1,2.$ Let $\tau_i={\rm ad}\, u_i\circ \tau,$ $i=1,2.$
Suppose that there is a sequence of unitaries $w_n\in U_0(M(B)/B$
such that $\lim_{n\to\infty} {\rm ad}\, w_n\circ
\tau_1(a)=\tau_2(a)\rforal a\in A. $ It follows that
$\lim_{n\to\infty} {\rm ad}\,u_2^*\circ {\rm ad}\,w_n\circ
\tau_1(a)=\tau(a)$ for all $a\in A.$ It follows from \ref{IIIP1} that
there is a subsequence $\{k(n)\}$ such that $\{[u_2^* w_{k(n)}
u_1]\}\in H_1^{ap}(K_0(A),  K_0(B)).$ This implies that
$\overline{[u_1]}=\overline{[u_2]},$ or $x_1=x_2.$
\end{proof}

\begin{Lem}\label{IIILKK}
Let $A$ be a separable unital \CA\, in ${\cal N}$ and $B$ be a
unital  purely infinite simple \CA.  Let $\tau: A\to
q_{\infty}(B)=l^{\infty}(B)/c_0(B)$ be a unital monomorphism and $\zeta\in H_1(K_0(A),
K_1(q_{\infty}(B))).$ Then, there is full monomorphism $\tau_1:
A\otimes C(S^1)\to q_{\infty}(B)$ such that
$$
[\tau_1|_{A}]=[\tau]\,\,\, {\rm in}\,\,\,
 KL(A, q_{\infty}(B))\andeqn \Gamma([\tau_1])([1_A\otimes
\imath)]=\zeta\, {\rm \,in} \,\,\,\,\,\,K_1(q_{\infty}(B)),
$$
where $\imath(z)=z$ for $z\in S^1.$
\end{Lem}

\begin{proof}
It follows from 8.5 and 7.7 of \cite{Lnef} that there is a
(group) isomorphism from the approximately unitary equivalence
classes of full monomorphisms from  $A$ to  $q_{\infty}(B)$ and
$KL(A, q_{\infty}(B)).$ Exactly the same proof of \ref{L1} proves
the theorem. One may also use the identification $KK^1(A,
Sq_{\infty}(B)))=KK(A, q_{\infty}(B))$ and apply \ref{L1}.
\end{proof}

\begin{thm}\label{Ttrm2}
Let $B$ be a $\sigma$-unital simple \CA\, with continuous scale
and  let $A$ be a unital separable amenable \CA\, in ${\cal N}.$

{\rm (i)} Let $\tau: A\to M(B)/B$ be a weakly
unital essential extension and $\{u_n\}\in U(M(B)/B)$
such that $[u_n]\in H_1^{ap}(K_0(A), K_1(M(B)/B))$
and
$$
\sigma(a)=\lim_{n\to\infty}{\rm ad}\, u_n\circ \tau(a)\rforal a\in A,
$$
where $\sigma: A\to M(B)/B$ is another essential extension. Then
$\sigma$ and $\tau$ are strongly unitarily equivalent.

{\rm (ii)} If $\prod_{n\in N}K_1(M(B)/B)=H_1^{ap}(K_0(A),
K_1(M(B)/B))$ (in particular, if \\
$K_1(M(B)/B)=H_1(K_0(A),
K_1(M(B)/B))$), then two weakly unital essential extensions
$\tau_1, \tau_2: A\to M(B)/B$ are  strongly approximately unitarily equivalent
if and only if $[\tau_1]=[\tau_2]$ in $KL(A, M(B)/B).$

{\rm (iii)} If $\tau: A\to M(B)/B$ is a weakly unital essential
extension, then there is an injection from\\
$K_1(M(B)/B)/H_1(K_0(A),K_1(M(B)/B)))$ to the set of strong
approximate unitary equivalence classes of weakly unital essential
extensions $\sigma$ for which $[\sigma]=[\tau]$ in $KL(A, M(B)/B)$
;

{\rm (iv)} In {\rm (iii)}, if furthermore,
$K_1(M(B)/B)/H_1(K_0(A), K_1(M(B)/B))$ is finite, then there is a
bijection from $K_1(M(B)/B)/H_1(K_0(A), K_1(M(B)/B))$ onto the set
of strong approximate unitary equivalence classes of weakly unital
essential extensions $\sigma$ for which $[\sigma]=[\tau]$ in
$KL(A, M(B)/B).$

{\rm (v)}  If neither of $\tau$ and $\sigma$ are weakly unital,
then $\tau$ and $\sigma$ are strongly approximately unitarily
equivalent if and only if $[\sigma]=[\tau]$ in $KL(A, M(B)/B).$

\end{thm}

\begin{proof}
Put $Q=M(B)/B.$ Then $Q$ is purely infinite and simple (see
\ref{PR}). In particular, $Q$ has $K_1$-$r$-cancellation (for
$r(n)=n$) and ${\rm cel}(M(B)/B)\le 2\pi +d$ (for any $d>0$) (see
\cite{P}).
 To see (i), let $\xi=\pi_{*1}([\{u_n\}]),$ where
$\pi: l^{\infty}(Q)\to q_{\infty}(Q).$ By passing to a
subsequence, without loss of generality, we may assume that
$\xi\in H_1(K_0(A), K_1(q_{\infty}(Q))).$ Let $\Phi: A\to
l^{\infty}(Q)$ be defined by $\Phi(a)=\{\sigma(a),\tau(a),\cdots,
\}$ for all $a\in A.$ Define ${\bar \tau}=\pi\circ \Phi.$ It
follows from \ref{IIILKK} that there exists a full monomorphism
${\bar \tau}_1: A\otimes C(S^1)\to q_{\infty}(Q)$ such that
$[({\bar \tau}_1)|_{A}]=[\Phi]$ in $KL(A,q_{\infty}(Q))$ and
$[{\bar \tau}_1( 1\otimes \imath)]=\xi.$ It follows from 7.7 of
\cite{Lneq} that $\Phi$ is approximately unitarily equivalent to
$({\bar \tau}_1)|_A.$ As in the proof of \ref{T1ext} we may assume
that $({\bar \tau}_1)|_A$ is unital. There is a sequence of
unitaries $w_n\in q_{\infty}(Q)$ such that
$$
\lim_{n\to\infty}{\rm ad}\, w_n\circ {\bar \tau}_1(a\otimes 1)=
\Phi(a)\rforal a\in A.
$$
Put $z={\bar\tau}_1(1\otimes \imath).$ Then
$$
\lim_{n\to\infty}\|w_n^*zw_n\Phi(a)-\Phi(a)w_n^*zw_n\|=0\rforal
a\in A.
$$
It is well known that there are unitaries $w_n(k), z(k)\in Q$ such
that $\pi(\{w_n(k)\})=w_n$ and $\pi(\{z(k)\})=z.$ Thus there is a
subsequence $\{n(k)\}$ such that
$$
\lim_{n\to\infty}\|w_{n(k)}(k)^*z(k)w_{n(k)}(k)\sigma(a)-\sigma(a)w_{n(k)}(k)^*z(k)w_{n(k)}\|=0
\rforal a\in A.
$$
Let $v_k=w_{n(k)}(k)^*z(k)w_{n(k)}(k).$ Then $[v_k]=[z(k)]$ in
$K_1(Q).$  By the definition of $\xi,$ one checks that
$[v_k]=[u_{n(k)}]$ for all sufficiently large $k.$ Note that
$v_k^*u_{n(k)}\in U_0(Q),$ since $Q$ is purely infinite and
simple. We have
$$
\lim_{n\to\infty} {\rm ad}\, u_{n(k)}v_k^*\circ
\tau(a)=\sigma(a)\rforal a\in A.
$$
This proves (i).

Note that (ii) follows from (i) and \ref{PT2} immediately.

(iii) follows from the same proof of \ref{IITappf}.

To see (iv), we note that, by (ii), we only need to show that the
injection is actually  surjective in this case. Suppose that
$\{u_n\}$ be a sequence of unitaries such that
$$
\sigma(a)=\lim_{n\to\infty}{\rm ad}\,u_n\circ \tau(a)\rforal a\in A.
$$
Let $\phi: K_1(M(B)/B)\to K_1(M(B)/B)/({\bar G}+H_1(K_0(A),
K_1(M(B)/B)))$ be the quotient map. By the assumption infinitely
many $\phi([u_n])$ are the same. Suppose that $\{k(n)\}$ is a
subsequence of $\N$ such that $\phi([u_{k(n)}])=x$ for $n=1,2,...$
for some $x\in K_1(M(B)/B)/({\bar G}+H_1(K_0(A), K_1(M(B)/B))).$
Choose $u\in U(M(B)/B$ such that $\phi([u])=x.$ Then
$$
\sigma(a)=\lim_{n\to\infty} {\rm ad}\,u_{k(n)}u^*\circ {\rm ad}\,
u\circ \tau(a)\rforal a\in A.
$$
Since  $[u_{k(n)}u^*]\in  H_1(K_0(A), K_1(M(B)/B)),$ by (i),
$\sigma$ is strongly unitarily equivalent to ${\rm ad}\, u\circ
\tau.$ This proves (iv).

For (v), by \ref{PT2}, we only need to show the ``if" part. Since
$M(B)/B$ is purely infinite simple \CA\, (see \cite{Lncs2}), there
is a unitary $v_1\in U(M(B)/B)$ such that
$v_1^*\tau_2(1)v_1=\tau_1(1).$ Put $e=1-\tau_1(1).$ Again, since
$M(B)/B$ is purely infinite and simple, there is a unitary
$v_2'\in e(M(B)/B)e$ such that $[v_2']=[v_1^*].$ Put
$v_2=v_2'+(1-e).$ Then $v_2^*v_1^*\tau_2(1)v_1v_2=\tau_1(1).$ Note
that $v_1v_2\in U_0(M(B)/B).$ Thus we may assume that
$\tau_1(1)=\tau_2(1).$ Let $u_n\in M(B)/B$ be a sequence of
unitaries such that
$$
\lim_{n\to\infty} {\rm ad}\, u_n\circ \tau_1(a)=\tau_2(a)\rforal
a\in A.
$$
Since $M(B)/B$ is purely infinite and simple, we obtain unitaries
$w_n\in e(M(B)/B)e$ such that $[w_n]=[u_n].$ Put
$z_n=w_n^*+(1-e).$ Then $z_nu_n\in U_0(M(B)/B).$ One verifies that
$$
\lim_{n\to\infty}{\rm ad}\, z_n\circ \tau_1(a)=\tau_2(a)\rforal
a\in A.
$$
\end{proof}

\begin{Remark}\label{RF}
{\rm  From the above theorem we also know that, for each $[\tau]$
in $KL(A, M(B)/B),$ there are at most $\prod_{n\in {\mathbb N}}
K_1(M(B)/B)/H^{ap}_1(K_0(A), K_1(M(B)/B))$ many different strong
approximate unitary equivalence classes of weakly unital essential
extensions $\sigma$ for which $[\sigma]=[\tau].$ However, given a
sequence of elements $\{x_n\}$  in $K_1(M(B)/B),$ we do not know
if there is a sequence of unitaries $\{u_n\}$ such that
$[u_n]=x_n$ and ${\rm ad}\, u_n\circ \tau(a)$ converges for any
$a\in A.$ This prevents us from determining  exactly how many
different strong approximate unitary equivalence classes of weakly
unital essential extensions $\sigma$ with $[\sigma]=[\tau].$ In
the case that $B=C\otimes {\cal K},$ one can make another
estimate. There are at most $|K_0(B)/H_1(K_0(A), K_0(B))|\cdot
|ext_{\mathbb Z}(K_*(A), K_*(B))|$ many strong approximate unitary
equivalence classes of weakly unital essential extensions $\sigma$
with $[\sigma]=[\tau].$ However, even if $[\sigma]\not=[\tau]$ in
$KK^1(A, B)$ but $[\sigma]=[\tau]$ in $KL^1(A, B),$ they could
still be strongly approximately  unitarily equivalent. For
example, this happens when $K_0(B)=H_1(K_0(A), K_0(B)),$ since in
this case, by \ref{IIITapp=}, all approximately unitarily
equivalent weakly unital full essential extensions are strongly
approximately unitarily equivalent. Full extensions in different
$KK$-classes could be strongly approximately unitarily equivalent
even in the case that $K_1(M(B)/B)/H_1(K_0(A),
K_1(M(B)/B))\not=\{0\}.$ For example, suppose that
$K_1(M(B)/B)/H_1(K_0(A), K_1(M(B)/B))$ is finite but $ext_{\mathbb
Z}(K_1(A), K_0(M(B)/B))$ is infinite. In this case there are only
finitely many strong approximate unitary equivalence classes of
weakly unital essential extensions which give the same element in
$KL^1(A, M(B)/B).$  Therefore there are infinitely many weakly
unital essential extensions in different $KK(A, M(B)/B))$ (but in
the same $KL(A, M(B)/B)$ ) that  are strongly approximately
unitarily equivalent.

On the other hand, the above results show that there are weakly
unital full essential extensions which give the same element in
$KK(A, M(B)/B)$  (or in $KK^1(A, B)$) may not strongly
approximately unitarily equivalent as long as
$K_1(M(B)/B)/H_1(K_0(A), K_1(M(B)/B))$ is not trivial. For
short,\\
$KK(A, M(B)/B)/KL(A,M(B)/B))$ can not be used to distinguish
strong approximate unitary equivalence from approximate unitary
equivalence. It is the group $H_1(K_0(A), K_1(M(B)/B)$ (or the
approximate version of it) that detects the difference.

}

\end{Remark}


\begin{thebibliography}{BK1}
{\small

\bibitem{B} B. Blackadar, {\em $K$-theory for Operator Algebras},
2nd ed. Mathematical Sciences Research Institute Publications, 5.
  Cambridge University Press, Cambridge, 1998.
\bibitem{Br2} L. G. Brown, {\em The universal coefficient theorem
for $Ext$ and quasidiagonality}, Operator Algebras and Group Representations, vol. 17,
Pitman Press, Boston, 1983, pp. 60-64.
\vspace{-0.08in}

\bibitem{BDF1}L. G. Brown, R. G. Douglas and P. A. Fillmore,
{\em Unitary equivalence modulo the compact operators and
extensions of $C\sp{*} $-algebras},
  Proceedings of a Conference on Operator Theory
(Dalhousie Univ., Halifax, N.S., 1973), pp. 58--128.
Lecture Notes in Math., Vol. 345, Springer, Berlin, 1973.
\vspace{-0.08in}

\bibitem{BDF3} L. G. Brown, R. G. Douglas and P. A. Fillmore,
{\em Extensions of $C\sp*$-algebras and $K$-homology},
 Ann. of Math. {\bf 105} (1977),
  265--324.
\bibitem{EK} G. A. Elliott and D. Kucerovsky,
{\em  An abstract Voiculescu-Brown-Douglas-Fillmore
absorption theorem},  Pacific J. Math.  {\bf 198}  (2001),
385--409.
\bibitem{GL} G. Gong and H. Lin, {\em Almost multiplicative morphisms
and $K$-theory},  Internat. J. Math.  {\bf 11}  (2000),   983--1000.
\bibitem{Lncs1} H. Lin, {\em Simple $C\sp *$-algebras with
continuous scales and simple corona algebras},
 Proc. Amer. Math. Soc.  112  (1991),
no. 3, 871--880.

\bibitem{Ln3} H. Lin, {\em $C^*$-algebra Extensions of $C(X)$},
           Memoirs Amer. Math. Soc.,  {\bf 115} (1995), no. 550.
\vspace{-0.08in}



\bibitem{LneII} H. Lin, {\em Extensions by $C^*$-algebras with real
rank zero II},
Proc. London
Math. Soc., 71 (1995), 641-674.
\vspace{-0.08in}


\bibitem{Lnamj} H. Lin, {\em Extensions of $C(X)$ by
simple $C^*$-algebras of real rank zero},
Amer. J. Math. {\bf 119} (1997), 1263-1289.
\vspace{-0.08in}

\bibitem{Lnb} H. Lin, {\em An Introduction to the  Classification
of Amenable $C^*$-Algebras}, World Scientific, 2001.
\vspace{-0.08in}

\bibitem{Lnpro} H. Lin, {\em A separable
Brown-Douglas-Fillmore Theorem and weak stability}, Trans. Amer.
Math. Soc., to appear.


\bibitem{Lncs2} H. Lin, {\em
Simple corona \CA s}, Proc. Amer. Math. Soc., to appear.

\bibitem{Lneq} H. Lin, {\em
Extensions by simple \CA s--quasidiagonal extensions}, Canad. J. Math.,
to appear.

\bibitem{Lnef} H. Lin, {\em Full extensions and approximate unitary equivalnce},
preprint, ArXiv. OA/0401242.

\bibitem{P} G. K. Pedersen, {\em $AW\sp *$-algebras and corona
$C\sp *$-algebras, contributions to noncommutative topology},  J. Operator
  Theory 15 (1986), 15--32.

\bibitem{PPV1} M. Pimsner, S. Popa and D.  Voiculescu,
{\em Homogeneous $C\sp{*} $-extensions of $C(X)\otimes K(H)$. I}, J.
  Operator Theory {\bf 1} (1979),  55--108.
\bibitem{Pc} N. C. Phillips, {\em A survey of exponential rank.
$C\sp *$-algebras: 1943--1993} (San Antonio, TX, 1993),  352--399,
Contemp. Math., 167, Amer. Math. Soc., Providence, RI, 1994.

\bibitem{RS} J.  Rosenberg and C.  Schochet, {\em The Kunneth
theorem and the universal coefficient theorem for Kasparov's
  generalized $K$-functor},  Duke Math. J. {\bf 55} (1987),
 431--474.



}
\end{thebibliography}
\end{document}